\documentclass[11pt]{article}
\newcommand		{\comment}[1]		{}

\usepackage{amsmath,amsthm} 			
\usepackage{mathtools}

\usepackage[left=2.54cm,right=2.54cm,top=2.54cm,bottom=2.54cm]{geometry}
\usepackage[nouppercase]{scrlayer-scrpage}
\usepackage{xspace}
\usepackage{enumitem}

\usepackage{avant} 				
\usepackage[sc,osf]{mathpazo} 			
	\AtBeginDocument{
		\DeclareSymbolFont{AMSb}{U}{msb}{m}{n}
		\DeclareSymbolFontAlphabet{\mathbb}{AMSb}
	}
\usepackage{extpfeil} 
\usepackage{mathabx}  
\usepackage{mathrsfs}	
\usepackage{fancyvrb}	

\usepackage{rotating,pdflscape}
\usepackage{afterpage}
\usepackage{graphicx}
\usepackage{float}
\usepackage{subcaption}
\usepackage{caption}
\newcommand{\mockalph}[1]{\!}

\usepackage[utf8]{inputenc}
\usepackage[T1]{fontenc}

\usepackage{tocloft}
\makeatletter
\renewcommand{\l@figure}{\@dottedtocline{1}{1em}{3.5em}}
\renewcommand{\l@table}{\@dottedtocline{2}{1em}{3.5em}}
\makeatother
\newcommand*{\noaddvspace}{\renewcommand*{\addvspace}[1]{}}
\addtocontents{lof}{\protect\noaddvspace}

\usepackage[hyphens]{url}		%
\usepackage{hyperref}	
\usepackage{cleveref} 	
\newcommand		{\myred}		{BrickRed}
\hypersetup{
			colorlinks		= true, 	
			urlcolor		= \myred, 
			linkcolor		= \myred, 	
			citecolor		= PineGreen 	
}
\newcommand		{\hyref}[1]		{\hyperref[#1]{\ref*{#1}}}

\newif\ifdebug
\ifdebug\usepackage{lineno}\linenumbers\else\fi

\usepackage{etoolbox}
\PassOptionsToPackage{dvipsnames,svgnames,x11names,table}{xcolor}

\newcommand		{\defm}[1]	{\textcolor{RoyalBlue}{#1}}

\usepackage[all]{xypic}\xyoption{rotate}
	\newdir{ >}{{}*!/-7pt/\dir{>}}
	\newdir{i}{{}*!/-7pt/\dir{(}	}
\usepackage{tikz}
\usetikzlibrary{matrix,fit,backgrounds,calc,arrows} 
\tikzstyle{image}=[rectangle,fill=Red!20,inner sep=-2pt]
\tikzstyle{nonzero}=[rectangle,fill=Navy!20,inner sep=0pt]
\tikzstyle{nonzerosm}=[rectangle,fill=Navy!20,inner sep=-2pt]

\usepackage{thmtools}
\makeatletter
\let\c@figure\c@table
\let\c@equation\c@table
\makeatother

\numberwithin{table}{section}
\numberwithin{figure}{section}
\newtheorem{abcthm}[table]{Theorem}

\newtheorem{theorem}[table]{Theorem}

\newtheorem{proposition}[table]{Proposition}

\newtheorem{corollary}[table]{Corollary}

\newtheorem{lemma}[table]{Lemma}

\newtheorem{claim}[table]{Claim}

\theoremstyle{definition}

\newtheorem{definition}[table]{Definition}

\newtheorem{notation}[table]{Notation}

\newtheorem{observation}[table]{Observation}

\newtheorem{conjecture}[table]{Conjecture}

\theoremstyle{remark}
\newtheorem{fact}[table]{Fact}

\newtheorem{example}[table]{Example}

\newtheorem{exercise}[table]{Exercise}

\newtheorem{problem}[table]{Problem}
\newtheorem{histrmks}[table]{Historical remarks}
\newtheorem{remark}[table]{Remark}
\newtheorem{remarks}[table]{Remarks}

\theoremstyle{plain}
\newtheorem*{thm*}{Theorem}
\newtheorem*{theorem*}{Theorem}
\newtheorem*{prop*}{Proposition}
\newtheorem*{proposition*}{Proposition}
\newtheorem*{lemma*}{Lemma}
\newtheorem*{corollary*}{Corollary}
\newtheorem*{cor*}{Corollary}

\theoremstyle{definition}
\newtheorem*{definition*}{Definition}
\newtheorem*{defn*}{Definition}
\newtheorem*{QQ*}{Question}
\newtheorem*{obs*}{Observation}
\newtheorem*{notation*}{Notation}
\newtheorem*{discussion*}{Discussion}

\theoremstyle{remark}
\newtheorem*{rmk*}{Remark}
\newtheorem*{remark*}{Remark}
\newtheorem*{examples*}{Examples}
\newtheorem*{example*}{Example}
\newtheorem*{EG*}{Example}
\newtheorem*{EGs*}{Examples}
\newtheorem*{fact*}{Fact}
\newtheorem*{prob*}{Problem}

\newcommand{\bthm}{\begin{theorem}}
\newcommand{\ethm}{\end{theorem}}
\newcommand{\bprop}{\begin{proposition}}
\newcommand{\eprop}{\end{proposition}}
\newcommand{\bcor}{\begin{corollary}}
\newcommand{\ecor}{\end{corollary}}
\newcommand{\bconj}{\begin{conjecture}}
\newcommand{\econj}{\end{conjecture}}
\newcommand{\blem}{\begin{lemma}}
\newcommand{\elem}{\end{lemma}}
\newcommand{\bclm}{\begin{claim}}
\newcommand{\eclm}{\end{claim}}
\newcommand{\bpf}{\begin{proof}}
\newcommand{\epf}{\end{proof}}
\newcommand{\bdetails}{\begin{details}}
\newcommand{\edetails}{\end{details}}
\newcommand{\bdefi}{\begin{definition}}
\newcommand{\edefi}{\end{definition}}
\newcommand{\bdefn}{\begin{definition}}
\newcommand{\edefn}{\end{definition}}
\newcommand{\bex}{\begin{example}}
\newcommand{\eex}{\end{example}}
\newcommand{\bprob}{\begin{problem}}
\newcommand{\eprob}{\end{problem}}
\newcommand{\bob}{\begin{observation}}
\newcommand{\eob}{\end{observation}}
\newcommand{\bexer}{\begin{exercise}}
\newcommand{\eexer}{\end{exercise}}
\newcommand{\bexers}{\begin{exercises}}
\newcommand{\eexers}{\end{exercises}}
\newcommand{\brmk}{\begin{remark}}
\newcommand{\ermk}{\end{remark}}
\newcommand{\bhist}{\begin{histrmks}}
\newcommand{\ehist}{\end{histrmks}}
\newcommand{\brmks}{\begin{remarks}}
\newcommand{\ermks}{\end{remarks}}
\newcommand{\bverb}{\begin{verbatim}}
\newcommand{\everb}{\end{verbatim}}
\newcommand{\bntn}{\begin{notation}}
\newcommand{\entn}{\end{notation}}
\newcommand{\bfct}{\begin{fact}}
\newcommand{\efct}{\end{fact}}
\newcommand{\bfcts}{\begin{facts}}
\newcommand{\efcts}{\end{facts}}
\newcommand{\benum}{\begin{enumerate}}
\newcommand{\eenum}{\end{enumerate}}
\newcommand{\bitem}{\begin{itemize}}
\newcommand{\eitem}{\end{itemize}}

\tracingpatches
\makeatletter
\patchcmd{\@setref}{\bfseries ??}{\bfseries\color{red} FIX ME!}{}{}
\patchcmd{\@setcite}{\bfseries ?}{\bfseries\color{red} FIX ME!}{}{}
\patchcmd{\@setcref}         {??}{\color{red} FIX ME!}{}{}
\patchcmd{\@setcref}         {??}{\color{red} FIX ME!}{}{}
\patchcmd{\@setcrefrange}    {??}{\color{red} FIX ME!}{}{}
\patchcmd{\@setcrefrange}    {??}{\color{red} FIX ME!}{}{}
\patchcmd{\@setcrefrange}    {??}{\color{red} FIX ME!}{}{}
\patchcmd{\@setcrefrange}    {??}{\color{red} FIX ME!}{}{}
\patchcmd{\@setcrefrange}    {??}{\color{red} FIX ME!}{}{}
\patchcmd{\@setcrefrange}    {??}{\color{red} FIX ME!}{}{}
\patchcmd{\@setnamecref}     {??}{\color{red} FIX ME!}{}{}
\patchcmd{\@setnamecref}     {??}{\color{red} FIX ME!}{}{}
\patchcmd{\@setcpageref}     {??}{\color{red} FIX ME!}{}{}
\patchcmd{\@setcpageref}     {??}{\color{red} FIX ME!}{}{}
\patchcmd{\@setcpagerefrange}{??}{\color{red} FIX ME!}{}{}
\patchcmd{\@setcpagerefrange}{??}{\color{red} FIX ME!}{}{}
\patchcmd{\@setcpagerefrange}{??}{\color{red} FIX ME!}{}{}
\patchcmd{\@setcpagerefrange}{??}{\color{red} FIX ME!}{}{}
\patchcmd{\@setcpagerefrange}{??}{\color{red} FIX ME!}{}{}
\patchcmd{\@cref}            {??}{\color{red} FIX ME!}{}{}
\def\blx@citation@entry#1#2{%
  \blx@bibreq{#1}%
  \ifinlist{#1}{\blx@cites}
    {}
    {\listgadd{\blx@cites}{#1}%
     \blx@auxwrite\@mainaux{}{\string\abx@aux@cite{#1}}}%
  \ifinlistcs{#1}{blx@segm@\the\c@refsection @\the\c@refsegment}
    {}
    {\listcsgadd{blx@segm@\the\c@refsection @\the\c@refsegment}{#1}}%
  \blx@ifdata{#1}%
    {}%
    {\ifcsdef{blx@miss@\the\c@refsection}%
       {\ifinlistcs{#1}{blx@miss@\the\c@refsection}%
          {{\bfseries\color{red} cite:} }%
          {\blx@logreq@active{#2{#1}}}}%
       {\blx@logreq@active{#2{#1}}}}}
\def\blx@citeadd#1{%
  \ifcsdef{blx@keyalias@\the\c@refsection @#1}
    {\edef\blx@realkey{\csuse{blx@keyalias@\the\c@refsection @#1}}}
    {\def\blx@realkey{#1}}%
  \expandafter\blx@citation\expandafter{\blx@realkey}\blx@msg@cundefon
  \expandafter\blx@ifdata\expandafter{\blx@realkey}
    {\advance\blx@tempcnta\@ne
     \listeadd\blx@tempa{\blx@realkey}}
    {\ifnum\blx@tempcntb>\z@\multicitedelim\fi
     \expandafter\abx@missing\expandafter{\blx@realkey}%
     \advance\blx@tempcntb\@ne}}
\makeatother

\newcommand{\presectionskip}{-1.5\baselineskip}
\newcommand{\postsectionskip}{0.3\baselineskip}
\makeatletter
\usepackage{titlesec}
\renewcommand{\section}{\@startsection
  {chapter}{0}{0mm}
  {\presectionskip}
  {\postsectionskip}
  {\sffamily\huge}}
\renewcommand{\section}{\@startsection
  {section}{1}{0mm}
  {\presectionskip}
  {\postsectionskip}
  {\sffamily\LARGE}}
\renewcommand{\subsection}{\@startsection
  {subsection}{2}{0mm}
  {\presectionskip}
  {\postsectionskip}
  {\sffamily\Large}}
\renewcommand{\subsubsection}{\@startsection
  {subsubsection}{3}{0mm}
  {\presectionskip}
  {\postsectionskip}
  {\sffamily\normalsize}}
\renewcommand{\@seccntformat}[1]{\csname the#1\endcsname.\quad}
\newcommand\HUGE{\@setfontsize\Huge{30}{47}} 
\makeatother
\usepackage{titling}
  \titleformat{\chapter}[display]
  {\sffamily\Large}
  {Chapter {\HUGE\normalfont\thechapter}}    
  {1em}
  {\huge}
  \titleformat{\section}
  {\sffamily\Large}
  {\normalfont{\@setfontsize\Large{18}{47}\thesection.}}    
  {1em}
  {}
  \titleformat{\subsection}
  {\sffamily\large}
  {\normalfont{\Large\thesubsection}.}    
  {1em}
  {}
  \titleformat{\subsubsection}
  {\sffamily\normalsize}
  {\normalfont{\large\thesubsubsection}.}    
  {1em}
  {}

\def			\SPSB#1#2			{\rlap{\textsuperscript{#1}}\textsubscript{#2}}
\makeatletter
\def			\smallunderbrace#1		{\mathop{\vtop{\m@th\ialign{##\crcr
							   $\hfil\displaystyle{#1}\hfil$\crcr
							   \noalign{\kern3\p@\nointerlineskip}%
							   \tiny\upbracefill\crcr\noalign{\kern3\p@}}}}\limits}
\makeatother

\makeatletter
\newcommand{\subalign}[1]{%
  \vcenter{%
    \Let@ \restore@math@cr \default@tag
    \baselineskip\fontdimen10 \scriptfont\tw@
    \advance\baselineskip\fontdimen12 \scriptfont\tw@
    \lineskip\thr@@\fontdimen8 \scriptfont\thr@@
    \lineskiplimit\lineskip
    \ialign{\hfil$\m@th\scriptstyle##$&$\m@th\scriptstyle{}##$\crcr
      #1\crcr
    }%
  }
}
\makeatother

\makeatletter
\newcommand		{\oset}[3][0ex]			{%
								\raisebox{.175ex}{$%
								  \mathrel{\mathop{#3}\limits^{
								    \vbox to#1{\kern-2\ex@
								    \hbox{$\scriptstyle#2$}\vss}}}
								    $}%
							    }
\makeatother

\usepackage{faktor,xfrac}

\newcommand		{\nd}		{\noindent}

\newcommand		{\mn}		{\mspace{-2mu}}

\newcommand		{\os}			{\overset}

\newcommand		{\mr}			{\mathrm}
\newcommand		{\bb}			{\mathbb}

\renewcommand		{\a}		{\alpha}

\renewcommand		{\b}		{\beta}

\newcommand		{\g}		{\gamma}

\renewcommand		{\epsilon}	{\varepsilon}

\renewcommand		{\l}		{\lambda}

\newcommand		{\s}		{\sigma}

\DeclareSymbolFont{cmletters}{OT1}{cmr}{m}{n}
\DeclareMathSymbol{\Ups}{\mathalpha}{cmletters}{"7}
\renewcommand		{\Upsilon}	{\Ups}

\newcommand		{\ceq}		{\coloneqq}

\let\union\cup%
\renewcommand		{\cup}		{\mspace{-1mu}\smile\mspace{-1mu}}
\makeatletter
\DeclarePairedDelimiterX
			{\pmodx}[1]	{(}{)}{{\operator@font mod}\mkern6mu#1}
					\renewcommand{\pmod}{%
					  \allowbreak
					  \if@display\mkern18mu\else\mkern8mu\fi
						  \pmodx
					}
\makeatother

\renewcommand		{\:}		{\colon}
\renewcommand		{\-}		{^{-1}}

\renewcommand		{\o}		{\circ}

\ExplSyntaxOn
\NewDocumentEnvironment{adjunctions}{O{}}
{
	\cs_set_eq:cN {@arraycr} \farin_arraycr:
	\keys_set:nn { farin/adjunction } { #1 }
	\begin{array}
		{
			@{ \hspace { \dim_eval:n { \l_farin_left_shift_dim + \l_farin_padding_dim } } }
			r
			@{ {\farin_strut:} \l_farin_symbol_tl {} }
			l
			@{ \hspace { \dim_eval:n { \l_farin_right_shift_dim + \l_farin_padding_dim } } }
		}
	}
	{
	\end{array}
}
\keys_define:nn { farin/adjunction }
{
	leftshift       .dim_set:N = \l_farin_left_shift_dim,
	leftshift       .initial:n = 0pt,
	rightshift      .dim_set:N = \l_farin_right_shift_dim,
	rightshift      .initial:n = 0pt,
	padding         .dim_set:N = \l_farin_padding_dim,
	padding         .initial:n = 6pt,
	symbol          .tl_set:N  = \l_farin_symbol_tl,
	symbol          .initial:n = \longrightarrow,
	verticalspacing .dim_set:N  = \l_farin_vertspac_dim,
	verticalspacing .initial:n = {3pt},
}
\cs_new_protected:Npn \farin_strut:
{
	\vrule height \dim_eval:n { \ht\strutbox + 1.2\l_farin_vertspac_dim }
	depth  \dim_eval:n { \dp\strutbox + \l_farin_vertspac_dim }
	width 0pt
}
\makeatletter
\exp_args:NNo \cs_new:Npn \farin_arraycr:
{
	\@arraycr\hline
}
\makeatother
\ExplSyntaxOff

\renewcommand		{\.}		{\cdot}

\newcommand		{\x}		{\times}

\DeclareMathOperator*	{\otimesvariable}{%
			\mathchoice {\raisebox{.85pt}{$\displaystyle\otimes$}}
						{\raisebox{.85pt}{$\otimes$}}
						{\raisebox{0.7pt}{$\scriptstyle\otimes$}}
						{\raisebox{0.2pt}{$\scriptscriptstyle\otimes$}}
						}
\newcommand		{\tensor}	{\otimesvariable}

\newcommand		{\ox}		{\tensor}

\newcommand		{\colim}	{\varinjlim}

\newcommand		{\exterior}	{\Lambda}
\newcommand		{\ext}		{\exterior}

\DeclareMathOperator	{\rk}		{rk }

\DeclareMathOperator	{\Tor}		{Tor}

\DeclareMathOperator	{\Ad}		{Ad }

\DeclareMathOperator	{\Aut}		{Aut }

\newcommand		{\U}		{\mr{U}}

\makeatletter
\newbox\xrat@below
\newbox\xrat@above
\newcommand		{\xrightarrowtail}[2][]	{%
						  \setbox\xrat@below=\hbox{\ensuremath{\scriptstyle #1}}%
						  \setbox\xrat@above=\hbox{\ensuremath{\scriptstyle #2}}%
				  \pgfmathsetlengthmacro{\xrat@len}{max(\wd\xrat@below,\wd\xrat@above)+.6em}%
  						\mathrel{\tikz [>->,baseline=-.55ex]
              					   \draw (0,0) -- node[below=-2pt] {\box\xrat@below}
                            					    node[above=-2pt] {\box\xrat@above}
                    						   (\xrat@len,0) ;}
						}
\newbox\xrat@below
\newbox\xrat@above
\renewcommand		{\xtwoheadrightarrow}[2][]{%
						  \setbox\xrat@below=\hbox{\ensuremath{\scriptstyle #1}}%
						  \setbox\xrat@above=\hbox{\ensuremath{\scriptstyle #2}}%
				  \pgfmathsetlengthmacro{\xrat@len}{max(\wd\xrat@below,\wd\xrat@above)+.6em}%
						 \mathrel{\tikz [->>,baseline=-.55ex]
					                 \draw (0,0) -- node[below=-2pt] {\box\xrat@below}
					                                node[above=-2pt] {\box\xrat@above}
						                       (\xrat@len,0) ;}
		       				}
\makeatother
\newcommand		{\xmono}	{\xrightarrowtail}
\newcommand		{\mono}		{\xmono{\phantom{\ \, }}}

\newcommand		{\xepi}		{\xtwoheadrightarrow}
\newcommand		{\epi}		{\xepi{\phantom{\ \, }}}

\newcommand		{\longto} 	{\longrightarrow}
\newcommand		{\lt}		{\longto}

\newcommand		{\lmt}		{\longmapsto}

\newcommand		{\longsimto}	{\os\sim\longto}
\newcommand		{\isoto}	{\longsimto}

\newcommand		{\iso}		{\cong}

\newcommand		{\Z}		{\bb Z}

\newcommand		{\C}		{\bb C}

\newcommand		{\K}		{K^*}

\newcommand		{\Gab}		{G^{\mathrm{ab}}}

\newcommand		{\Gad}		{G^{\Ad}}

\newcommand		{\KG}		{\K_G}

\let\lim\varprojlim

\newcommand{\bAd}{\b^{\Ad}}
\newcommand		{\Gbi}{G^{\mathrm{bi}}}

\begin{document}
	\title{The K-theory of the conjugation action}
	\author{Jeffrey D.~Carlson}
	\maketitle

\begin{abstract}
  In 1999, Brylinski and Zhang computed the complex equivariant
  K-theory of the conjugation self-action of a compact, connected Lie group
  with torsion-free fundamental group. In this note we show it is possible to
  do so in under a page.
\end{abstract}

Brylinski and Zhang~\cite{brylinskizhang} showed that if $G$ is a compact,
connected Lie group with torsion-free fundamental group, then the equivariant
K-theory of its conjugation action $\defm \Gad$ is isomorphic to the ring
$\defm{\Omega^*_{RG/\Z}}$ of Grothendieck differentials on the complex
representation ring $\defm{RG}$ of $G$. Their proof uses results on holomorphic
differentials on complex manifolds, a reduction to the case $G$ is a torus, and
some algebraic geometry. We show a more concrete and arguably more natural
expression for the ring $\KG(\Gad)$ can be obtained rapidly using only
Hodgkin's K\"unneth spectral sequence~\cite{hodgkin1975kunneth}, in the same
manner they already use it, and elementary algebraic considerations. We then
show this purely algebraic isomorphism admits a satisfying geometric
interpretation, and remark finally that this geometric version gives back
Brylinski and Zhang's description in terms of Grothendieck differentials at no
added cost.

\begin{theorem*}[Brylinski--Zhang {\cite[Theorem~3.2]{brylinskizhang}}]
Let $G$ be a compact, connected Lie group with torsion-free fundamental group.
Then $\KG(\Gad)$
is isomorphic to $RG \ox \K G$ as an $RG$-algebra.
Under this identification,
the forgetful map $\defm f\mn\: \KG(\Gad) \lt \K G$
becomes reduction with respect to the augmentation ideal~$\defm{IG}$ of $RG$.
\end{theorem*}
\begin{proof}
  Write $\defm\Gbi$ for $G$ under the $(G\x G)$-action $(h,k)\.g =
  hgk\-$. The orbit space of $\Gbi \x \Gbi$ under the restricted, free diagonal
  action of $1 \x G$ is $(G \x 1)$-equivariantly diffeomorphic to $\Gad$ via
  $(g',g) \lmt g'g\-$, so when $ X = Y = \Gbi$, Hodgkin's $(\Z \x \Z/2)$-graded
  K\"{u}nneth spectral sequence $\Tor_{R(G \x G)}(\K_{G \x G} X, \K_{G \x G} Y)
  \implies \K_{G \x G} (X \x Y)$ reduces to $\Tor_{RG \ox RG}(RG,RG) \implies
  \K_G(\Gad)$. Here the two structure maps $RG \ox RG \lt RG$ are the both the
  multiplication of $RG$. Recall~\cite[Prop.~11.1]{hodgkin1975kunneth} that
  under our hypotheses, $RG$ is the tensor product of a polynomial ring on
  generators $\defm{y_i} \in IG$ and a Laurent polynomial ring on generators
  $\defm{t_j} \in 1 + IG$. Let $\defm P$ be the free abelian group on
  generators $\defm{q_i}$ and $\defm{w_j}$ and let $\defm\g\: P \lt IG$ be the
  linear map taking $q_i$ to $y_i$ and $w_j$ to $t_j - 1$. Then an $(RG)^{\ox
  2}$-module resolution of $RG$ is given by $RG \ox \ext P \ox RG$, with
  differential the derivation vanishing on $RG \ox \Z \ox RG$ and sending $1
  \ox z \ox 1$, for $z \in P$, to~$1 \ox 1 \ox \g(z) - \g(z) \ox 1 \ox 1$. To
  compute the Tor, apply $- \ox_{RG \ox RG} RG$ to this resolution to obtain
  the \textsc{cdga} $RG \ox \ext P$ with $0$ differential, $RG$ in bidegree
  $(0,0)$, and $P$ in bidegree $(-1,0)$.

The spectral sequence collapses because the differentials $d_r$ for $r \geq 2$
send all generators into the right half-plane. Since $\pi_1 G$ is torsion-free
and $X = G$ is locally contractible of finite covering dimension, the spectral
sequence strongly converges to the intended target. Hence $RG \ox \ext P$ is
the graded algebra associated to a filtration $(\defm{F_p})_{p \leq 0}$ of
$\KG(\Gad)$ with $F_0 \iso RG$ and $F_{-1}/F_0 \iso RG \ox P$. Since $RG$ and
$\ext P$ are both free abelian, there is no additive extension problem, so
$\KG(\Gad)$ is also free abelian as a group. Let $\widetilde z_k$ be elements
in $F_{-1}$ lifting $1 \ox q_i$ and $1 \ox w_j$ under the
isomorphism~$F_{-1}/F_0 \iso RG \ox P$. Then the $\widetilde z_k$ anticommute
with each other because they lie in $K^1_G(\Gad)$ and square to $0$ since
$\KG(\Gad)$ contains no $2$-torsion, and by induction, they generate
$\KG(\Gad)$ as an $RG$-algebra, so $\KG(\Gad) = RG \ox \ext[\widetilde z_k]$.

To see the forgetful map $f$ is as claimed, note that forgetting the $(G \x
1)$-action on $\Gbi$ induces a map to the spectral sequence $\Tor_{RG}(\Z,\Z)
\implies \K G$, which again collapses by lacunary considerations. Computing
$\Tor_{RG}(\Z,\Z) \iso \ext P$ with the resolution $\ext P \ox RG $ of $\Z$
shows the map $E_2(f)\colon RG \otimes \ext P \longrightarrow \ext P$ is
reduction modulo $IG$ and $\ext P \iso \K G$.  
\end{proof}
\begin{remark*}
  We can be completely explicit about the exterior generators. As
 Hodgkin~\cite[Thm.~A]{hodgkin1967lie}
  observed, the injection $\U(n) \mono
  \U \ceq \colim \U(n)$ induces an additive map $\defm\b\: RG \lt K^1 G$
  descending to a group isomorphism between the module $IG/(IG)^2$ of
  indecomposables of $RG$ and the module $P\K G$ of primitives of the exterior
  Hopf algebra $\K G \iso \ext P\K G$. In particular, a set of generators is
  given by $\b(\l_i) = \b(\l_i - \dim \l_i)$ for $\defm{\l_i}$ lifts in $G$ of
  the fundamental representations of the commutator subgroup $G'$ and $\b(t_j)
  = \b(t_j - 1)$ for $\defm{t_j}\: G \to G/G' \epi \U(1)$ circular coordinate
  functions of the torus $\Gab = G/G' \iso (S^1)^{\rk G - \rk G'}$. Let $\defm
  Q = \{\l_i,t_j\}_{i,j}$.

The map $\b$ in fact factors as $f \o \bAd$ for a map $\defm{\bAd}\: RG \lt
\KG(\Gad)$, already giving surjectivity of $f$ since $P \K G$ generates $\K G =
\ext P \K G$ as a ring. Atiyah~\cite[Lem.~2, pf.]{atiyah1965lie} described
$\bAd$ and hence $\b$ geometrically: given a representation $\rho\: G \lt
\U(n)$, we can build a representative $\defm E$ of $\b(\rho)$ via the clutching
construction, taking with two trivial bundles $CG \x \C^n$ over the cone
$\defm{CG}$ on $G$ and gluing them along $G \x \C^n$ via the relation $(g,v)
\sim \big(g,\rho(g)v\big)$ to obtain a bundle over the suspension $CG \union_G
CG$. The action $h\.(g,v) = \big(hgh^{-1},\rho(h)v\big)$ of $G$ on $\Gad \x
\C^n$ preserves this relation and so induces a $G$-action on $E$ making it a
$G$-equivariant bundle over the suspension of~$\Gad$.

The $RG$-module structure on $\KG(X)$ is always given by $\s\.[E] = \big[\a(\s)
\ox E\big]$, where if $\s\: G \lt \Aut V$ is a representation, then
$\defm{\a}(\s)$ is the trivial bundle $X \x V$ equipped with the diagonal
$G$-action. Thus $\sigma \ox \prod_k \b(\rho_k) \lmt
\a(\s)\.\prod_k\bAd(\rho_k)$ for $\rho_k \in Q$ gives an explicit isomorphism
$RG \ox \K G \isoto \KG(\Gad)$.  
\end{remark*}

\smallskip

\begin{remark*}
  The first paragraph of our proof is a variant of
  Brylinski--Zhang's \S4~\cite{brylinskizhang}. Once the ring structure is
  determined as in our proof's second paragraph, replacing \S\S5--6, one knows
  their map~$\smash{ \defm\phi\: {\Omega^*_{RG/\Z}} \lt \KG(\Gad)}$ from the
  ring of Grothendieck differentials is an isomorphism as soon as one knows it
  is a well-defined $RG$-algebra map~\cite[Prop~3.1]{brylinskizhang}, for the
  class $\bAd(\rho) \in K^1_G(\Gad)$ from the previous remark is in fact the
  same as\footnote{\ the corrected version of---see
  Fok~{\cite[Rmk.~2.8.1]{fok2014real}}} Brylinski--Zhang's $\phi(d\rho)$, so
  that $f \o \phi$ takes a basis~$\{\defm{d\rho} : \rho \in Q\}$ of the free
  $RG$-module $\Omega^1_{RG/\Z}$ to a $\Z$-basis $\big\{\b(\rho):\rho \in
  Q\big\}$ of $P\K G$.
\end{remark*}

{\footnotesize\bibliography{bibshort} }

\nd\footnotesize{%
	\url{jeffrey.carlson@tufts.edu
	}
}
\end{document}